\documentclass[10pt]{amsart}

\usepackage{amsmath}
\usepackage{amsthm}
\usepackage{amsfonts}
\usepackage{amssymb,latexsym}
\usepackage[all]{xy}
\usepackage{graphicx}
\usepackage[mathscr]{eucal}
\usepackage{verbatim}

\theoremstyle{plain}
    \newtheorem{thm}{Theorem}[section]
    
    \newtheorem{lemma}[thm]{Lemma}

\theoremstyle{definition}

\theoremstyle{remark}
    \newtheorem{rem}[thm]{Remark}
    \newtheorem{example}[thm]{Example}

\newtheorem*{namedtheorem}{\theoremname}
\newcommand{\theoremname}{testing}

\newcommand{\A}{\mathcal A}
\newcommand{\B}{\mathcal B}
\newcommand{\C}{\mathcal C}

\newcommand{\rar}{\ensuremath{\rightarrow}}

\newcommand{\n}{\noindent}
\newcommand{\length}{\text{length}}
\newcommand{\MaxSpec}{\text{MaxSpec}}

\DeclareMathOperator{\Spec}{Spec}

\begin{document}

%
\title{Classifying subcategories of modules over a PID.}
\date{\today}

\author{Sunil K. Chebolu}
\address {Department of Mathematics \\
          University of Western Ontario \\
          London, ON. N6A 5B7 }
\email{schebolu@uwo.ca}


\begin{abstract} Let $R$ be a commutative ring. A full additive subcategory $\C$ of $R$-modules is
\emph{triangulated} if whenever two terms of a short exact sequence
belong to $\C$, then so does the third term. In this note we give a
classification of triangulated subcategories of finitely generated
modules over a principal ideal domain. As a corollary we show that
in the category of finitely generated modules over a PID, thick
subcategories (triangulated subcategories closed under direct
summands), wide subcategories (abelian subcategories closed under
extensions) and Serre subcategories (wide subcategories closed under
kernels) coincide and correspond to specialisation closed subsets of
$\Spec(R)$.
\end{abstract}

\maketitle

\section{Introduction}

A full additive subcategory $\C$ of $R$-modules  is said to be an
\emph{exact} subcategory if it is closed under extensions of
$R$-modules. These categories were introduced by Quillen in
\cite{quillen-ktheory1} where he defined them more generally in an
abstract and intrinsic manner. However, in all examples of interest,
exact categories embed in some abelian category and so it is best to
think of them as defined above. It is clear from the work of Quillen
\cite{quillen-ktheory1, quillen-ktheory2}, Waldhausen \cite{wald}
and other recent model category theorists \cite{dugg-ship} that
exact categories play an important role in algebraic $K$-theory and
Morita theory. So a good understanding of them can be very helpful.

In this paper we consider some  enriched versions  of exact
subcategories: triangulated subcategories (exact subcategories that
are closed under kernels of surjective maps and cokernels of
injective maps), thick subcategories (triangulated subcategories
that are closed under direct summands), wide subcategories (abelian
subcategories that are closed under extensions), Serre classes (wide
subcategories that are closed under submodules and quotient
modules), and torsion theories (Serre classes that are closed under
arbitrary direct sums) are some examples of subcategories that have
been studied. So one has the following hierarchy of exact
subcategories.
\begin{center}
Torsion theory $\Rightarrow$ Serre class $\Rightarrow$  Wide
$\Rightarrow$ Thick $\Rightarrow$ Triangulated.
\end{center}
(See \cite[section 1]{wide} for a good account on these
subcategories.) Triangulated subcategories are quite general. Note,
for instance, that a triangulated subcategory need not be closed
under kernels, cokernels or direct summands. The following simple
example should illustrate this point.

\begin{example} The subcategory $\C$ of all even dimensional rational vector spaces in the category of
finite dimensional rational vector spaces is readily seen to be a
triangulated subcategory. However, it is not closed under kernels.
Consider the projection map
\[
\begin{array}{ccc}
\mathbb{Q}^2 & \longrightarrow &  \mathbb{Q}^2 \\
(x,y) &  \longmapsto  & (x,0).
\end{array}
\]
This map has a one dimensional kernel and a one dimensional
cokernel, and hence the subcategory $\C$ is neither closed under
kernels nor cokernels. It is also clear that $\C$ is not closed
under direct summands.
\end{example}

The terms "triangulated subcategory" and "thick subcategory" are
widely used in the context of triangulated categories such as
derived categories of rings. We have defined analogues of these
subcategories in category of modules and have decided to give them
the same names. However, this should not cause any confusion to the
reader because in this paper we work exclusively in the category of
$R$-modules.

 Since triangulated categories form the most general family in
the hierarchy of exact subcategories,  if we have a classification
of triangulated subcategories, we can use that classification to
derive classifications of other families of subcategories. However,
classifying triangulated subcategories of modules over a general
ring seems to be hopelessly difficult. Even when restricted to
finitely generated modules, we do not know of a recipe for
classifying triangulated subcategories analogous to the K-theory
recipe of Thomason \cite{Th} for classifying the triangulated
subcategories of a triangulated category. Nevertheless, when $R$ is
a principal ideal domain (PID), in view of the structure theorem for
finitely generated $R$-modules there is a hope for such a
classification. And as we will see, one can give a complete
classification of all the aforementioned subcategories over PIDs.
Finally we should point out that the material in this paper should
be in the reach of any undergraduate student who is familiar with
basic properties of modules over a PID.

The paper is organised as follows. We begin in the next section with
a quick recap of some basic facts about modules over PIDs and the
notion of length of a module. We then classify the triangulated
subcategories of finitely generated modules over a PID in section 3.
In the last section we derive classifications of the aforementioned
families of subcategories over PIDs.

\vskip 4mm \noindent \textbf{Acknowledgement:}  I would like  to
thank John Palmieri and Sverre Smal{\o} for some interesting
discussions on this problem.

\section{Principal ideal domains}
In this short section we recall  some well-known facts about PIDs.
We begin with the structure theorem for finitely generated (f.g.)
modules. Let $R$ be a PID. If $M$ is a finitely generated
$R$-module, then
\[ M \cong  F(M)  \oplus \left( \bigoplus_{p\, \in \,\Spec(R)}  E_p(M) \right),\]
where $F$ is a free $R$-module of finite rank and $E_p(M)$ denotes
the $R$-submodule of $M$ which consists of all elements that are
annihilated by some power of $p$.
Further, for each prime element $p$ that occurs in the above sum,
$E_p(M)$ decomposes as
\[E_p(M) \cong R/(p^{l_1}) \oplus R/(p^{l_2}) \oplus \cdots \oplus R/(p^{l_s}),\]
with uniquely determined sequence $1 \leq l_1 \leq l_2 \leq \cdots \leq l_s $.

Now we recall the notion of the length of a module. For an $R$-module $M$, a chain
\[M=M_0 \supsetneq M_1 \supsetneq M_2 \supsetneq \cdots \supsetneq M_r = 0\]
is called a composition series if each $M_i/M_{i+1}$ is a simple
module (one that does not have any non-trivial submodules). The
length $l(M)$ of a module $M$ that is both Artinian and Noetherian
is defined to be the length of any composition series of $M$. The
fact that this is well-defined is part of the Jordan-Holder theorem.
If $0 \rightarrow N \rightarrow M \rightarrow K \rightarrow 0 $ is a
short exact sequence of $R$-modules of finite length, then we have
$l(M) = l(N) + l(K)$. In other words, the length function $l(\;)$ is
an Euler characteristic function on the category of modules of
finite length.
Finally by rank of a finitely generated module we will mean the rank
of its free part, and similarly by $p$-length of a module, we will
mean the $p$-length of its $p$-torsion part.

\section{Classification of triangulated subcategories over a PID.}
Henceforth $R$ will denote a principal ideal domain unless otherwise
stated.  We now classify all the triangulated subcategories in the
category $\A$ of finitely generated modules over  $R$. We first set
up some notations.  We will denote by $\B$ the subcategory of
finitely generated torsion $R$-modules. For every prime ideal $p$ in
$R$, define Euler characteristic functions $\chi_0$ and  $\chi_p$ on
$\A$ and $\B$ respectively as follows:
\[
\chi_p(X) :=
 \left\{
   \begin{array}{ll}
      \dim_{R_{0}}( X \otimes_R R_{0})  & \mbox{ when $p =0$}, \\
      \length (X \otimes_R R_{(p)}) & \mbox{ when $p \ne 0$}
    \end{array}
 \right.
\]

These Euler characteristic functions define some full subcategories
of $\A$ and $\B$ as follows. For each positive integer $k$, define
$I_k$ as
\[ I_k := \{ X \in \A: \chi_0(X) \equiv 0 \mod k \}, \]
and for each subgroup $H$ of ${\underset{p \, \in \,
\MaxSpec(R)}{\oplus}}\mathbb{Z}$ that is  generated by elements all
of whose components are non-negative, define $J_H$ as
\[J_H := \{ X \in \B: {\underset{p\, \in \, \MaxSpec(R)}{\oplus}} \; \chi_p(X) \in H\}. \]
Since $\chi_p$ is an Euler characteristic function, it is clear that
all these subcategories are triangulated subcategories. We now show
that these are all the triangulated subcategories in $\A$. We begin
with two lemmas which will help streamline the proof.

\begin{lemma}\emph{(Descending lemma)} Let $\C$ be a triangulated
subcategory of $\A$ and let  $M$ be a module in $\C$ with $p$-length
$r$. Then there exists another module $N$ in $\C$ whose $p$-torsion
is $(R/p)^r$ and whose $q$-torsion is identical with that of $M$ for
all primes $q \ne p$.
\end{lemma}

\begin{proof}
Let $M \cong \oplus_{i=1}^{k} R/p^{r_i} \oplus L$, with $\Sigma \,
r_i = r$ and $L$ $p$-torsion free. We want to generate modules with
lower highest $p$-order (In the module $M$, the highest $p$-order is
$max\{r_i\}$.) by keeping the $p$-length of the module constant.
This is done until we get a module with highest $p$-order $=1$, or
equivalently, until we get $\oplus_{i=1}^{r}
R/p \oplus L$. For better clarity, we break the proof  into two steps. \\

\n (a) Assume without loss of generality that $r_1 = max\{r_i\}$. We
first  build $R/p \oplus R/p^{r_1-1}$ from $R/p^{r_1}$. The
following pair of short exact sequences will do the job. (For
clarity the obvious quotient maps are not labelled.)
\[ 0 \rar R/p^{r_1} \overset{(1,p)}{\rar} R/p^{r_1-1} \oplus R/p^{r_1+1} \rar R/p^{r_1} \rar 0\]
\[ 0 \rar R/p^{r_1} \overset{(0,p)}{\rar} R/p^{r_1-1} \oplus R/p^{r_1+1} \rar R/p^{r_1-1} \oplus R/p  \rar 0\]

\n (b) Now we build $R/p \oplus R/p^{r_1-1} \oplus R/p^{r_2} \oplus
\cdots \oplus R/p^{r_k} \oplus L $ from $M = R/p^{r_1} \oplus
R/p^{r_2} \oplus \cdots \oplus R/p^{r_k} \oplus L$. (Note that this
decreases the highest $p$-order by one.) This can be obtained easily
by adding  a split short exact sequence
\[ 0 \rar G \rar G \oplus G \rar G \rar 0\]
with  $G= \oplus_{i=2}^{k} R/p^i \oplus L$ to the pair of short
exact sequences in part (a):
\[ 0 \rar R/p^{r_1} \oplus G \rar R/p^{r_1-1} \oplus R/p^{r_1+1} \oplus \; G \oplus \; G \rar R/p^{r_1} \oplus G \rar 0\]
\[ 0 \rar R/p^{r_1} \oplus G \rar R/p^{r_1-1} \oplus R/p^{r_1+1} \oplus\; G \oplus  \; G\rar R/p^{r_1-1} \oplus R/p \oplus G \rar 0\]

 A straightforward downward induction will decrease the highest
order to $1$. In other words, it produces $\oplus_{i=1}^{l} R/p
\oplus L$, completing the proof of the lemma.

\end{proof}

\begin{lemma}\emph{(Ascending lemma)} Let $\C$ be a triangulated
subcategory of $\A$, and let $M$ be a module in $\C$ whose
$p$-torsion part is $(R/p)^r$. Now given any partition $\Sigma \,
r_i$ of $r$ into non-negative integers, then there exists another
module $N$ in $\C$ whose $p$-torsion is $\oplus \, R/p^{r_i}$ and
whose $q$-torsion is identical with that of $M$ for all primes $q
\ne p$.
\end{lemma}

\begin{proof}
We have to start with $(R/p)^r \oplus L$ where $L$ is $p$-torsion
free and construct the module $\oplus (R/p^{r_i}) \oplus L$, where
the exponents $r_i$ are such that  $\Sigma r_i = r$. Again for
clarity, we break the construction into two steps.

\n (a) We first generate $R/p^{r+1}$ from $R/p^r \oplus R/p$. The
following pair of short exact sequences will do this job.
\[ 0 \rar (R/p \oplus R/p^r) \overset{h}{\rar} (R/p \oplus  R/p^{r+1} )\oplus R/p^{r} \rar R/p \oplus R/p^r \rar 0\]
\[ 0 \rar (R/p \oplus R/p^r) \hookrightarrow (R/p \oplus R/p^{r}) \oplus R/p^{r+1} \rar R/p^{r+1} \rar 0\]
The map $h$ sends $(x ,y)$ to $(x, py, 0)$. The remaining maps are
the obvious inclusion and quotient maps.

 \n (b) Now we have to show that $R/p^{r+1} \oplus G$ can be
generated  from $R/p^r \oplus R/p \oplus G$ where $G$ is any
$R$-module. As before we add, to the above pair of short exact
sequences, a split short exact sequence $0 \rar G \rar G \oplus G
\rar G \rar 0$ to get
\[ 0 \rar R/p \oplus R/p^r \oplus G \rar R/p \oplus R/p^{r+1} \oplus R/p^{r} \oplus G \oplus G \rar R/p \oplus R/p^r \oplus G\rar 0\]
\[ 0 \rar R/p \oplus R/p^r \oplus G \rar R/p \oplus R/p^{r} \oplus R/p^{r+1} \oplus G \oplus G \rar R/p^{r+1} \oplus G\rar 0\]

Again a simple induction will complete the proof of the lemma.
\end{proof}

\begin{thm} Let $R$ be a PID and let $\A$ denote the category of finitely generated $R$-modules.
Then a non-zero subcategory $\C$ of $\A$ is triangulated  if and
only if either $\C = I_k$ for some positive integer $k$, or $J_H$
for some subgroup $H$ of $\underset{p \, \in \, \MaxSpec(R)}{\oplus}
\mathbb{Z}$ that is generated by elements all of whose components
are non-negative.
\end{thm}

\begin{proof} First of all it is easy to see that the subcategories
$I_k$ and $J_H$ as defined in the statement of the theorem are
triangulated and pair-wise distinct subcategories. So we have to
show that every triangulated subcategory $\C$ is equal to one of
these. The proof  of this
statement divides naturally into two cases.\\

\n \textbf{Case(i)} \emph{$\C$ contains a module of rank at least
one.} Pick a module $M$  of smallest non-zero rank (exists by
assumption) and let $k$ denote the rank of $M$. We claim that $\C =
I_k$. In other words, we have to show that $\C$ comprises of all
modules of the form $R^{kl} \oplus T$ where $l$ is a non-negative
integer and $T$ is a torsion module. This will be done by following
a series of straightforward reductions.  First of all, it suffices
to build free modules and torsion modules from $M$ (because
triangulated subcategories are closed under taking direct sums).
Then, the exact sequence
\[ 0 \rar Tor(M) \rar M(= R^k \oplus Tor(M)) \rar R^k \rar 0\]
gives a further reduction ($Tor(M)$ denote the torsion submodule of
$M$): It suffices to build an arbitrary torsion module out of $M$.
Now recall from the structure theorem that any torsion module is a
direct sum of cyclic modules of the form $R/p^t$. Again,
triangulated subcategories are closed under taking direct sums,
therefore it is enough to produce $R/p^t$ for any prime $p$ and any
integer $t$. The following short exact sequence tells us that this
is always possible: (here we use the fact that $k\ge1$)
\[ 0 \rar R \oplus( R^{k-1} \oplus Tor{M} ) \;\; {\overset{p^k \oplus (id \oplus id) }{\longrightarrow}} \;\; R \oplus (R^{k-1} \oplus Tor{M})
\rar R/p^t \rar 0.\]
So this completes the first case. \\

\n \textbf{Case(ii)}  \emph{All modules in $\C$ have rank zero.}
Equivalently, this means that $\C$ consists of torsion $R$-modules.
Let $\chi(-) := \oplus_{p \in \, \MaxSpec(R)} \; \chi_p(-)$ and
define $H$ to be the subgroup  of$\underset{p \, \in \,
\MaxSpec(R)}{\oplus} \mathbb{Z}$ generated by the elements $\chi(X)$
as $X$ ranges in $\C$. We claim that $\C = J_H$. So let $M$ be an
$R$-module with $\chi(M)$ belonging to $H$. We want to show that $M$
in $ \C$.  By the ascending lemma it suffices to show that there is
a module $N$ in $\C$ which is of the form $\oplus_{p \, \in \,
\MaxSpec(R)} (R/p)^{r_p}$ and has the same $p$-lengths as those of
$M$ for all primes $p$. By assumption, we know that
\[ \chi(M) = \sum a_i\; \chi(A_i) - \sum b_i \;\chi(B_i),\]
for some modules $A_i$ and $B_i$ in $\C$ and some positive
coefficients $a_i$ and $b_i$. By the descending lemma $A_i$ and
$B_i$ can be chosen such that
\[ A_i = \bigoplus_{p \, \in \, \MaxSpec(R)} (R/p)^{\alpha_{p}}, \]
\[ B_i = \bigoplus_{p \, \in \, \MaxSpec(R)} (R/p)^{\beta_{p}}. \]
where $\alpha_i$ and $\beta_i$ are non-negative integers. Note that
for each $p$, $\sum b_i \beta_{p} - a_i \alpha_{p}$ is  the
$p$-length of the module $M$ and therefore is a non-negative
integer. So there is an  inclusion $\oplus (B_i)^{b_i}
\hookrightarrow \oplus (A_i)^{a_i}$ whose cokernel is evidently the
desired candidate for $N$.
\end{proof}

\section{Classifications of other families of subcategories}

Recall the definitions of thick, wide and Serre subcategories from
the introduction. It turns out that these three families coincide in
$\A$, the category of finitely generated modules over a PID. This
will be shown  directly from the classification of triangulated
subcategories of $\A$ that is obtained in the previous section.

\subsection{Thick subcategories} To get a classification of the non-zero thick
subcategories of $\A$ we have to find out which of the $I_k$s and
the $J_H$s are closed under direct summands. One can easily check
that $I_k$ is closed under direct summands if and only if  $k=1$ and
$J_H$ is so precisely when $H$ is spanned by the unit vectors of the
form $e_p$ which have zeros everywhere except in the one spot
(corresponding to the prime $p$) where there is a one. In other
words, the only thick subcategories of $\A$ are $\A$ itself and the
subcategories of $S$-torsion modules where $S$ is a subset of
maximal primes in $R$.

\subsection{Wide subcategories} To get a classification of the non-zero wide
subcategories, we look for those thick subcategories that are closed
under kernels and cokernels. It is clear that both $\A$ and the
category of $S$-torsion modules (where $S$ is a subset of maximal
primes) are both abelian and closed under extensions and therefore
wide.

\subsection{Serre subcategories} It is easy to see that both  $\A$
and the category of $S$-torsion modules ($S$ is a subset of maximal
primes) are closed under subobjects and quotient objects. So these
are the only Serre subcategories in $\A$. \\

Thus we have show that in $\A$ the thick subcategories, wide
subcategories and Serre subcategories all coincide.

It is also clear that these subcategories are in 1-1 correspondence
with the specialisation closed subsets (subsets that are a union of
closed subsets in the Zariski topology) of  $\Spec(R)$ because a
non-trivial specialisation closed subset of  $\Spec(R)$ (when $R$
 is a PID) is precisely a subset of maximal primes,  and the trivial ones
 being $\Spec(R)$ and the empty set which correspond respectively to
 the categories $\A$ and the trivial category.
\begin{rem} There are some interesting classifications of some of these subcategories in the literature. For
example, the wide subcategories of finitely presented modules over a
regular coherent ring have been classified by Hovey \cite{wide}, and
the thick subcategories of finite dimensional representations over
finite $p$-groups have been classified by Benson-Carlson-Rickard
\cite{bcr}.
\end{rem}

\bibliographystyle{alpha}

\begin{thebibliography}{BCR97}

\bibitem[BCR97]{bcr}
D.~J. Benson, Jon~F. Carlson, and Jeremy Rickard.
\newblock Thick subcategories of the stable module category.
\newblock {\em Fund. Math.}, 153(1):59--80, 1997.

\bibitem[DS04]{dugg-ship}
Daniel Dugger and Brooke Shipley.
\newblock {$K$}-theory and derived equivalences.
\newblock {\em Duke Math. J.}, 124(3):587--617, 2004.

\bibitem[Gra76]{quillen-ktheory2}
Daniel Grayson.
\newblock Higher algebraic {$K$}-theory. {II} (after {D}aniel {Q}uillen).
\newblock In {\em Algebraic $K$-theory (Proc. Conf., Northwestern Univ.,
  Evanston, Ill., 1976)}, pages 217--240. Lecture Notes in Math., Vol. 551.
  Springer, Berlin, 1976.

\bibitem[Hov01]{wide}
Mark Hovey.
\newblock Classifying subcategories of modules.
\newblock {\em Trans. Amer. Math. Soc.}, 353(8):3181--3191 (electronic), 2001.

\bibitem[Qui73]{quillen-ktheory1}
Daniel Quillen.
\newblock Higher algebraic {$K$}-theory. {I}.
\newblock In {\em Algebraic $K$-theory, I: Higher $K$-theories (Proc. Conf.,
  Battelle Memorial Inst., Seattle, Wash., 1972)}, pages 85--147. Lecture Notes
  in Math., Vol. 341. Springer, Berlin, 1973.

\bibitem[Tho97]{Th}
R.~W. Thomason.
\newblock The classification of triangulated subcategories.
\newblock {\em Compositio Math.}, 105(1):1--27, 1997.

\bibitem[Wal85]{wald}
Friedhelm Waldhausen.
\newblock Algebraic {$K$}-theory of spaces.
\newblock In {\em Algebraic and geometric topology (New Brunswick, N.J.,
  1983)}, volume 1126 of {\em Lecture Notes in Math.}, pages 318--419.
  Springer, Berlin, 1985.

\end{thebibliography}

\end{document}